\documentclass[10pt,leqno]{article}

\usepackage{amsmath,amssymb,amsthm,mathrsfs,dsfont}

\usepackage[margin=3cm]{geometry} 

\usepackage{titlesec,hyperref}

\usepackage{color}

\usepackage{fancyhdr}
\titleformat{\subsection}{\it}{\thesubsection.\enspace}{1pt}{}

\newtheorem{theo}{Theorem}[section]
\newtheorem{lemm}[theo]{Lemma}
\newtheorem{defi}[theo]{Definition}
\newtheorem{coro}[theo]{Corollary}

\newtheorem{rema}[theo]{Remark}
\numberwithin{equation}{section}

\allowdisplaybreaks 

\newcommand\ep{{\varepsilon}} 
\newcommand\lm{{\lesssim}}

\begin{document}
\title{The Liouville theorem and the $L^2$ decay of the FENE dumbbell model of polymeric flows
\hspace{-4mm}
}

\author{Wei Luo$^1$
\quad Zhaoyang Yin$^2$ \\[10pt]
Department of Mathematics, Sun Yat-sen University,\\
510275, Guangzhou, China.\\[5pt]
}
\footnotetext[1]{Email: \it luowei23@mail2.sysu.edu.cn}
\footnotetext[2]{Email: \it mcsyzy@mail.sysu.com.cn}
\date{}
\maketitle
\hrule

\begin{abstract}
In this paper we mainly investigate the finite extensible nonlinear elastic (FENE) dumbbell model with dimension $d\geq2$ in the whole space. We first proved that there is only the trivial solution for the steady-state FENE model under some integrable condition. Our obtained results generalize and cover the classical results to the stationary Navier-Stokes equations. Then, we study about the $L^2$ decay of the co-rotation FENE model. Concretely, the $L^2$ decay rate of the velocity is $(1+t)^{-\frac{d}{4}}$ when $d\geq3$, and $\ln^{-k}{(e+t)}, k\in \mathds{N}^{+}$ when $d=2$. This result improves considerably the recent result of \cite{Schonbek2} by Schonbek. Moreover, the decay of general FENE model has been considered.\\

\vspace*{5pt}
\noindent {\it 2010 Mathematics Subject Classification}: 35Q30, 76B03, 76D05, 76D99.

\vspace*{5pt}
\noindent{\it Keywords}: The FENE dumbbell model; The Liouville theorem; Stationary solution; The Navier-Stokes equation; The $L^2$ decay.
\end{abstract}

\vspace*{10pt}

\tableofcontents

\section{Introduction}
   In this paper we consider the finite extensible nonlinear elastic (FENE) dumbbell model \cite{Bird}:
   \begin{align}
\left\{
\begin{array}{ll}
u_t+(u\cdot\nabla)u-\nu\Delta{u}+\nabla{P}=div~\tau, ~~~~~~~div~u=0,\\[1ex]
\psi_t+(u\cdot\nabla)\psi=div_{R}[- \sigma(u)\cdot{R}\psi+\beta\nabla_{R}\psi+\nabla_{R}\mathcal{U}\psi],  \\[1ex]
\tau_{ij}=\int_{B}(R_{i}\nabla_{j}\mathcal{U})\psi dR, \\[1ex]
u|_{t=0}=u_0,~~\psi|_{t=0}=\psi_0, \\[1ex]
(\beta\nabla_{R}\psi+\nabla_{R}\mathcal{U}\psi)\cdot{n}=0 ~~~~ \text{on} ~~~~ \partial B(0,R_{0}) .\\[1ex]
\end{array}
\right.
\end{align}
In (1.1)~~$\psi(t,x,R)$ denotes the distribution function for the internal configuration and $u(t,x)$ stands for the velocity of the polymeric liquid, where $x\in\mathbb{R}^{d}$ and $d\geq2$ means the dimension. Here the polymer elongation $R$ is bounded in ball $ B=B(0,R_{0})$ of $\mathbb{R}^{d}$ which means that the extensibility of the polymers is finite. $\beta=\frac{2k_BT_a}{\lambda}$, where $k_B$ is the Boltzmann constant, $T_a$ is the absolute temperature and $\lambda$ is the friction coefficient. $\nu>0$ is the viscosity of the fluid, $\tau$ is an additional stress tensor and $P$ is the pressure. The Reynolds number $Re=\frac{\gamma}{\nu}$ with $\gamma\in(0,1)$ and the density $\rho=\int_B\psi dR$. Moreover the potential $\mathcal{U}(R)=-k\log(1-(\frac{|R|}{|R_{0}|})^{2})$ for some constant $k>0$. $\sigma(u)$ is the drag term. In general, $\sigma(u)=\nabla u$. For the co-rotation case, $\sigma(u)=\nabla u-(\nabla u)^{T}$.

   This model describes the system coupling fluids and polymers. The system is of great interest in many branches of physics, chemistry, and biology, see \cite{Bird,Masmoudi.W}. In this model, a polymer is idealized as an "elastic dumbbell" consisting of two "beads" joined by a spring that can be modeled by a vector $R$. At the level of liquid, the system couples the Navier-Stokes equation for the fluid velocity with a Fokker-Planck equation describing the evolution of the polymer density. This is a micro-macro model (For more details, one can refer to  $\cite{Bird}$, $\cite{Doi}$, $\cite{Masmoudi.W}$ and $\cite{Masmoudi.G}$).

%
%
%


In the paper we will take $\beta=1$ and $R_{0}=1$.
Notice that $(u,\psi)$ with $u=0$ and $$\psi_{\infty}(R)=\frac{e^{-\mathcal{U}(R)}}{\int_{B}e^{-\mathcal{U}(R)}dR}=\frac{(1-|R|^2)^k}{\int_{B}(1-|R|^2)^kdR},$$
is a trivial solution of (1.1). By a simple calculation, we can rewrite (1.1) for the following system:
\begin{align}
\left\{
\begin{array}{ll}
u_t+(u\cdot\nabla)u-\nu\Delta u+\nabla{P}=div~\tau,  ~~~~~~~div~u=0,\\[1ex]
\psi_t+(u\cdot\nabla)\psi=div_{R}[-\sigma(u)\cdot{R}\psi+\psi_{\infty}\nabla_{R}\frac{\psi}{\psi_{\infty}}],  \\[1ex]
\tau_{ij}=\int_{B}(R_{i}\nabla_{R_j}\mathcal{U})\psi dR, \\[1ex]
u|_{t=0}=u_0, \psi|_{t=0}=\psi_0, \\[1ex]
\psi_{\infty}\nabla_{R}\frac{\psi}{\psi_{\infty}}\cdot{n}=0 ~~~~ \text{on} ~~~~ \partial B(0,1) .\\[1ex]
\end{array}
\right.
\end{align}
{\bf Remark.} As in the reference \cite{Masmoudi.G}, one can deduce that $\psi=0$ on the boundary.

Let us recall that the historical Liouville theorem (LT) states that a bounded entire holomorphic function is constant. This property can be generalized to a linear homogeneous elliptic system. However, whether the (LT) holds for a usual nonlinear elliptic system is hard to answer. The famous problem is the (LT) for stationary Navier-Stokes (SNS) equation. For $d=2$, (LT) for (SNS) was proved in
\cite{Gilbarg} by Gilbarg and Weinberger, while for $d=4$ is obtained by Galdi in
\cite{Galdi}. As far as we know, for $d=3$, this is still an open problem. The earliest result is due to Galdi \cite{Galdi} under the additional condition $u$ belongs to $L^{\frac{9}{2}}(\mathbb{R}^3)$. Recently, Chae and Yoneda \cite{Chae1} proved the (LT) for (SNS) if $u$ has a suitable behavior at infinity. In \cite{Chae2}, Chae obtained the result under the condition $u$ belongs to $W^{2,\frac{6}{5}}(\mathbb{R}^3)$. For the axially symmetric Navier-Stoke equation, Korobkov, Pileckas and Russo show the (LT) if the solutions are in absence of swirl.

To our best knowledge, there are no any results about the Liouville theorem for stationary FENE model (1.2). In this paper, we investigate the Liouville theorem for (1.2) with $d\geq2$. By using the similar idea as in \cite{Galdi} and \cite{Gilbarg}, we obtain the desire result for (1.2) under the responding integrable condition. If $d=3$, we add some additional condition which is different with that mentioned in \cite{Chae1}, \cite{Chae2} and \cite{Gilbarg}. Moreover, our result can be reduced to the Liouville theorem for Navier-Stoke equation and generalizes the result in \cite{Gilbarg}.

 In \cite{Schonbek1}, Schonbek proved the $L^2$ decay of the velocity for the Navier-Stoke equation and obtained the decay rate $(1+t)^{-\frac{d}{4}}$ which is in accord with that of the heat equation, this is a very interesting result. Recently, Schonbek \cite{Schonbek2} studied about the $L^2$ decay of the velocity for the co-rotation FENE dumbbell model, and obtained the
 decay rate $(1+t)^{-\frac{d}{4}+\frac{1}{2}}$. Moreover, she guessed that the correct decay rate should be $(1+t)^{-\frac{d}{4}}$ however she cannot use the bootstrap argument as in \cite{Schonbek1} because of the additional stress tensor. In this paper, we improved this result and verified that the $L^2$ decay rate is $(1+t)^{-\frac{d}{4}}$ with $d\geq 3$ i.e. Schonbek's guess is right. If $d=2$, Schonbek's result did not give the decay, and we proved that the decay rate is $\ln^{-k}(e+t)$ for any $k\geq 0$. The main idea is that we toke a parameter in the $L^2$ energy estimate such that the bootstrap argument is valid. Moreover, we also studied about the $L^2$ decay for the general FENE dumbbell model.

  The paper is organized as follows. In Section 2 we introduce some notations and  give some preliminaries which will be used in the sequel. In Section 3 we prove the Liouville theorem for the stationary FENE model. In Section 4 we study about the $L^2$ decay for FENE model by using the Fourier splitting method.

\section{Notations and preliminaries}
  In this section we first introduce some notations that we shall use throughout the paper.

For $p\geq1$, we denote by $\mathcal{L}^{p}$ the space
$$\mathcal{L}^{p}=\big\{\psi \big|\|\psi\|^{p}_{\mathcal{L}^{p}}=\int \psi_{\infty}|\frac{\psi}{\psi_{\infty}}|^{p}dR<\infty\big\}.$$

  We will use the notation $L^{p}_{x}(\mathcal{L}^{q})$ to denote $L^{p}[\mathbb{R}^{d};\mathcal{L}^{q}]:$
$$L^{p}_{x}(\mathcal{L}^{q})=\big\{\psi \big|\|\psi\|_{L^{p}_{x}(\mathcal{L}^{q})}=(\int_{\mathbb{R}^{d}}(\int_{B} \psi_{\infty}|\frac{\psi}{\psi_{\infty}}|^{q}dR)^{\frac{p}{q}}dx)^{\frac{1}{p}}<\infty\big\}.$$
When $p=q$, we also use the short notation $\mathcal{L}^p$ for $L^p_x(\mathcal{L}^{p})$ if there is no ambiguity.

The symbol $\widehat{f}=\mathcal{F}(f)$ denotes the Fourier transform of $f$.

Moreover, we denote by $\dot{\mathcal{H}}^1$ the space
$$\dot{\mathcal{H}}^1=\big\{g\big| \|g\|_{\dot{\mathcal{H}}^1}=(\int_B|\nabla_R g|^2\psi_\infty dR)^{\frac{1}{2}}\big\}.$$
Sometimes we write $f\lm g$ instead of $f\leq Cg$ where $C$ is a constant. We agree that $\nabla$ stands for $\nabla_x$ and $div$ stands for $div_x$.

 If the function spaces are over $\mathbb{R}^d$ and $B$ with respect to the variable $x$ and $R$, for simplicity, we drop $\mathbb{R}^d$ and $B$ in the notation of function spaces if there is no ambiguity.

The following lemma allows us to estimate the extra stress tensor $\tau$.
\begin{lemm}\cite{Masmoudi.G}\label{Lemma1}
There exists a constant $C$ such that for $\psi\geq0$ and $\sqrt{\frac{\psi}{\psi_\infty}}\in \dot{\mathcal{H}}^1$, we have
\begin{align}
|\tau|^2\leq C\bigg(\int_B\psi dR\bigg)\bigg(\int_B|\nabla_R\sqrt{\frac{\psi}{\psi_\infty}}|^2\psi_\infty dR\bigg).
\end{align}
\end{lemm}

\begin{lemm}\cite{Masmoudi.W}\label{Lemma2}
 If $\int_B \psi dR=0$ and
   $\displaystyle\int_B\bigg|\nabla _R  \bigg(\displaystyle\frac{\psi}{\psi _\infty }\bigg)\bigg|^2{\psi _\infty} dR<\infty$ with $p\geq 2$, then there exists a constant $C$ such that
   \[\int_{B}\frac{|\psi|^{2}}{\psi_{\infty}}dR\leq C \displaystyle\int_B\bigg|\nabla _R  \bigg(\displaystyle\frac{\psi}{\psi _\infty }\bigg)\bigg|^2{\psi _\infty} dR.\]
\end{lemm}
\begin{lemm}\label{Lemma3}
\cite{Masmoudi.W} For all $\varepsilon>0$, there exists a constant $C_{\varepsilon}$ such that
$$|\tau|^2\leq\varepsilon\int_{B}\psi_{\infty}|\nabla_{R}\frac{\psi}{\psi_{\infty}}|^{2}dR
+C_{\varepsilon}\int_{B}\frac{|\psi|^{2}}{\psi_{\infty}}dR.$$
\end{lemm}
\section{The Liouville theorem}
 In this section, we assume that $\sigma(u)=\nabla u$. Let us define the suitable stationary weak solution for (1.2).
\begin{defi}
 A couple of functions $(u,\psi)$ with $div~u=0$ is called a suitable stationary weak solution for (1.2) if the following conditions hold
 \begin{align}\tag{I}
 u\in [\dot{H}^1(\mathbb{R}^d)]^d,~~~\nabla_R\sqrt{\frac{\psi}{\psi_\infty}}\in L^2(\mathbb{R}^d\times B, \psi_\infty dxdR), ~~~\psi\geq 0,
 \end{align}
 \begin{align}\tag{II}
 \frac{\psi}{\psi_\infty}(\ln\frac{\psi}{\psi_\infty}-1)\in L^1(\mathbb{R}^d\times B, \psi_\infty dxdR),
 \end{align}
 \begin{align}\tag{III}
 \lim_{|x|\rightarrow\infty}u(x)=0,~~~~\lim_{|x|\rightarrow\infty}\psi(x,R)=\psi_{\infty},
 \end{align}
\begin{align}\tag{IV}
&\int_{\mathbb{R}^{d}}[(u\otimes u):\nabla v+P\cdot div~v] dx=\int_{\mathbb{R}^{d}}(\tau:\nabla v +\nu \nabla u:\nabla v) dx,~~\forall v\in C^\infty_0(\mathbb{R}^{d}), \\
\nonumber &\int_{\mathbb{R}^{d}\times B}u\psi\cdot\nabla_{x}\phi dxdR =\int_{\mathbb{R}^{d}\times B}[-\nabla{u}\cdot{R}\psi+\psi_{\infty}\nabla_{R}\frac{\psi}{\psi_{\infty}}]\cdot\nabla_{R}\phi dxdR,~~\forall \phi\in C^\infty_0(\mathbb{R}^d\times B).
\end{align}
\end{defi}
{\bf Remark.} The definition 3.1 is associated with the definition in \cite{Masmoudi.G} which corresponds to the evolution equations. The condition (I) is to ensure the regularity of the weak solution, while the condition (II) is called the entropy condition.

Our main results are the following:
\begin{theo}\label{th1}
Let $(u,\psi)$ be a bounded suitable stationary weak solution to (1.2) in $\mathbb{R}^d$. Assume that $\int_B\psi dR=1$ and there exist two constants $C_1,~C_2$ such that $0<C_1\leq\frac{\psi}{\psi_\infty}\leq C_2$. If
\begin{align}
u\in [L^{\frac{3d}{d-1}}(\mathbb{R}^d)]^d,
\end{align}
 then $u=0$ and $\psi=\psi_\infty$.
\end{theo}
\begin{theo}\label{th2}
Let $(u,\psi)$ be a bounded suitable stationary weak solution to (1.2) in $\mathbb{R}^3$. Assume that $\int_B\psi dR=1$ and there exist two constants $C_1,~C_2$ such that $0<C_1\leq\frac{\psi}{\psi_\infty}\leq C_2$. Let $1\leq p_i,~q_i,~r_i<\infty~ (i=1,2,3)$ and if
\begin{align}
u_i\in L^{p_i}_{x_1}L^{q_i}_{x_2}L^{r_i}_{x_3},~~ \text{with} ~~~\frac{1}{p_i}+\frac{1}{q_i}+\frac{1}{r_i}=\frac{2}{3},~~(i=1,2,3),
\end{align}
then $u=0$ and $\psi=\psi_\infty$.
\end{theo}
\begin{rema}
By the Sobolev embedding theorem, we have $\dot{H}^1(\mathbb{R}^d)\hookrightarrow L^{\frac{2d}{d-2}}(\mathbb{R}^d)$. If $d\geq 4$, one can see that $\frac{3d}{d-1}\geq \frac{2d}{d-2}$, which implies that $L^\infty(\mathbb{R}^d)\cap L^{\frac{2d}{d-2}}(\mathbb{R}^d)\hookrightarrow L^{\frac{3d}{d-1}}(\mathbb{R}^d)$. Hence, we can get rid of the condition (1.3) in Theorem 1.2  when $d\geq 4$.
\end{rema}
The proof of Theorem \ref{th2} is valid for Navier-Stokes equation and we have the following result:
\begin{coro}\label{co}
Let $u$ be a bounded stationary weak solution to the Navier-Stoke equation in $\mathbb{R}^3$. Assume that $u\in \dot{H}^1(\mathbb{R}^3)$ and $\lim_{|x|\rightarrow\infty}u(x)=0$. Let $1\leq p_i,~q_i,~r_i< \infty~ (i=1,2,3)$ and if
\begin{align}
u_i\in L^{p_i}_{x_1}L^{q_i}_{x_2}L^{r_i}_{x_3},~~ \text{with} ~~~\frac{1}{p_i}+\frac{1}{q_i}+\frac{1}{r_i}=\frac{2}{3},~~(i=1,2,3),
\end{align}
then $u=0$.
\end{coro}
\begin{rema}
By taking $p_i=q_i=r_i=\frac{9}{2}$ in Corollary \ref{co}. Our corollary cover the result in \cite{Galdi}. If we take $p_i=q_i=6$ which implies that $r_i=3$ i.e. $u\in L^6_{x_1}L^6_{x_2}L^3_{x_3}$. Since $\dot{H}^1(\mathbb{R}^3)\hookrightarrow L^6(\mathbb{R}^3)$, it follows that we only add the integrable condition $L^3$ in $x_3$. 
\end{rema}
\subsection{The proof of Theorem \ref{th1}}
In this subsection, we begin to prove Theorem \ref{th1}. For $K>0$, choose $\eta_K(x)$ to be a positive smooth cut-off function satisfying:
\begin{align}
\eta_K(x)=1,~~\text{if}~~|x|\leq K,~~\eta_K(x)=0,~~\text{if}~~|x|\geq 2K,~~|\nabla\eta_K(x)|\leq \frac{C}{K}, \text{for some constant}~C.
\end{align}
Since $u$ is bounded, it follows by density argument, for each fix $K>0$ we may choose $v_K(x)=u(x)\eta_K(x)$ as a test function, then we have
\begin{align}
\int_{\mathbb{R}^d}[(u\otimes u):\nabla (u\eta_K)]+P div(u\eta_K) dx =\int_{\mathbb{R}^d} \tau :\nabla (u \eta_K)+\nu\nabla u:\nabla (u\eta_K) dx.
\end{align}
Notice that $div~u=0$. By virtue of integration by parts, we compute that
\begin{align}
\nu\int_{\mathbb{R}^d}|\nabla u|^2\eta_K dx&=-\nu\int_{\mathbb{R}^d}\nabla u\cdot u\cdot\nabla\eta_K+\int_{\mathbb{R}^d}|u|^2u\cdot\nabla\eta_K dx+\frac{1}{2}\int_{\mathbb{R}}\nabla(|u|^2)\cdot u\eta_K dx+\int_{\mathbb{R}^d}P u\cdot \nabla \eta_K dx\\
\nonumber&-\int_{\mathbb{R}}\tau\cdot u\cdot \nabla\eta_K dx-\int_{\mathbb{R}}\tau:\nabla u \eta_K dx\\
\nonumber&=-\nu\int_{\mathbb{R}^d}\nabla u\cdot u\cdot\nabla\eta_K+\frac{1}{2}\int_{\mathbb{R}^d}|u|^2u\cdot\nabla\eta_K dx+\int_{\mathbb{R}^d}P u\cdot \nabla \eta_K dx\\
\nonumber&-\int_{\mathbb{R}^d}\tau\cdot u\cdot \nabla\eta_K dx-\int_{\mathbb{R}^d}\tau:\nabla u \eta_K dx.
\end{align}
Note that $\nabla\eta_K=0$ if $|x|\leq K$ and $|x|\geq 2K$. Then, we deduce from the above equality that
\begin{align}
\nu\int_{|x|\leq K}|\nabla u|^2dx&\leq \nu\int_{K\leq |x|\leq 2K}|\nabla u||u||\nabla \eta_K|dx+\frac{1}{2}\int_{K\leq |x|\leq 2K}|u|^3||\nabla \eta_K|dx+\int_{K\leq |x|\leq 2K}|\tau||u||\nabla \eta_K|dx\\
\nonumber&+\int_{K\leq |x|\leq 2K}|P||u||\nabla \eta_K|dx-\int_{\mathbb{R}^d}\tau:\nabla u \eta_K dx\\
\nonumber&=I^K_1+I^K_2+I^K_3+I^K_4-\int_{\mathbb{R}^d}\tau:\nabla u \eta_K dx.
\end{align}
Now we estimate the terms $I^K_1$ to $I^K_4$. If $d>2$, by virtue of H\"{o}lder's inequality with index $p=2$, $q=\frac{2d}{d-2}$ and  $r=d$, we have
\begin{align}\label{3.8}
I^K_1\lm \frac{1}{K}\int_{K\leq |x|\leq 2K}|\nabla u||u| dx\lm \frac{1}{K}\|\nabla u \|_{L^2(K\leq |x|\leq 2K)}\|u\|_{L^{\frac{2d}{d-2}}(K\leq |x|\leq 2K)}|K|\lm \|\nabla u \|^2_{L^2(K\leq |x|\leq 2K)}.
\end{align}
If $d=2$, by virtue of H\"{o}lder's inequality with index $p=2$ and $q=2$, we deduce that
\begin{align}
I^K_1\lm \frac{1}{K}\int_{K\leq |x|\leq 2K}|\nabla u||u| dx\lm \frac{1}{K}\|\nabla u \|_{L^2(K\leq |x|\leq 2K)}\|u\|_{L^{\infty}}|K|\lm \|u\|_{L^\infty}\|\nabla u \|_{L^2(K\leq |x|\leq 2K)}.
\end{align}
Using H\"{o}lder's inequality with index $p=\frac{d}{d-1}$, $q=d$, we obtain
\begin{align}
I^K_2\lm \frac{1}{K}\int_{K\leq |x|\leq 2K}|u|^3 dx\lm \|u\|^3_{L^{\frac{3d}{d-1}}(K\leq |x|\leq 2K)}.
\end{align}
By the same argument as $I^K_1$, we see that
\begin{align}
&I^K_3\lm  \|\tau\|_{L^2(K\leq |x|\leq 2K)}\|\nabla u\|_{L^2(K\leq |x|\leq 2K)}, ~~\text{if}~~~ d>2, \\
\nonumber &I^K_3\lm  \|\tau\|_{L^2(K\leq |x|\leq 2K)}\|u\|_{L^\infty}, ~~~~~~~~~~~~~~~~~\text{if}~~~ d=2.
\end{align}
Taking advantage of Lemma \ref{Lemma1} and using the fact that $\int_B\psi dR=1$, yield that
\begin{align}
\|\tau\|_{L^2(K\leq |x|\leq 2K)} \leq (\int_{\{K\leq |x|\leq 2K\}\times B} \bigg|\nabla_R\sqrt{\frac{\psi}{\psi_\infty}}\bigg|^2\psi_\infty dxdR)^{\frac{1}{2}}.
\end{align}
Since $P=\sum_{1\leq i,j\leq d}\mathfrak{R}^i\mathfrak{R}^j(u_iu_j-\tau_{ij})$ with $\mathfrak{R}$ is the usual Risez operator and using the fact that $\|\mathfrak{R}^if\|_{L^p}\leq \|f\|_{L^p}$, we have
\begin{align}\label{3.13}
&I^K_4\lm \|u\|^3_{L^{\frac{3d}{d-1}}(K\leq |x|\leq 2K)}+\|\tau\|_{L^2(K\leq |x|\leq 2K)}\|\nabla u\|_{L^2(K\leq |x|\leq 2K)}, ~~\text{if}~~~ d>2, \\
\nonumber &I^K_4\lm \|u\|^3_{L^{\frac{3d}{d-1}}(K\leq |x|\leq 2K)}+ \|\tau\|_{L^2(K\leq |x|\leq 2K)}\|u\|_{L^\infty}, ~~~~~~~~~~~~~~~~~\text{if}~~~ d=2.
\end{align}
From (\ref{3.8})-(\ref{3.13}), we deduce that $\lim_{K\rightarrow\infty}I^K_i=0,~~i=1,2,3,4$.

Thanks to $C_1\leq \frac{\psi}{\psi_\infty}\leq C_2$, for each $K$ we may choose $\phi_K(x)=\ln{\frac{\psi}{\psi_\infty}}\eta_K(x)$ as a test function to get
\begin{align}\label{3.14}
\int_{\mathbb{R}^{d}\times B}u\psi\cdot\nabla_{x}(\ln{\frac{\psi}{\psi_\infty}}\eta_K(x)) dxdR =\int_{\mathbb{R}^{d}\times B}[-\nabla{u}\cdot{R}\psi+\psi_{\infty}\nabla_{R}\frac{\psi}{\psi_{\infty}}]\cdot\nabla_{R}(\ln{\frac{\psi}{\psi_\infty}}\eta_K(x)) dxdR.
\end{align}
By directly calculating, we see that
\begin{align}\label{3.15}
\int_{\mathbb{R}^{d}\times B}\psi_\infty\nabla_{R}\frac{\psi}{\psi_{\infty}}\cdot\nabla_{R}(\ln{\frac{\psi}{\psi_\infty}}\eta_K(x))dxdR&=\int_{\mathbb{R}^{d}\times B}\psi_\infty\frac{\psi_\infty}{\psi}\bigg|\nabla_{R}\frac{\psi}{\psi_{\infty}}\bigg|^2\eta_K(x)dxdR\\
\nonumber&=4\int_{\mathbb{R}^{d}\times B}\psi_\infty\bigg|\nabla_{R}\sqrt{\frac{\psi}{\psi_{\infty}}}\bigg|^2\eta_K(x)dxdR.
\end{align}
Plugging (\ref{3.15}) into (\ref{3.14}) yields
\begin{multline}\label{3.16}
4\int_{\mathbb{R}^{d}\times B}\psi_\infty\bigg|\nabla_{R}\sqrt{\frac{\psi}{\psi_{\infty}}}\bigg|^2\eta_K(x)dxdR=\int_{\mathbb{R}^{d}\times B}u\psi\cdot\nabla_{x}(\ln{\frac{\psi}{\psi_\infty}}\eta_K(x)) dxdR\\
+\int_{\mathbb{R}^{d}\times B}\nabla{u}\cdot{R}\psi\nabla_{R}(\ln{\frac{\psi}{\psi_\infty}}\eta_K(x))dxdR.
\end{multline}
Since $\nabla_{x}\ln{\frac{\psi}{\psi_\infty}}=\nabla_{x}(\ln\psi-\ln\psi_\infty)=\frac{\nabla_{x}\psi}{\psi}$ and $div~u=0$, it follows that
\begin{align}\label{3.17}
\int_{\mathbb{R}^{d}\times B}u\psi\cdot\nabla_{x}(\ln{\frac{\psi}{\psi_\infty}}\eta_K(x)) dxdR&=\int_{\mathbb{R}^{d}\times B} (u\nabla_x \psi \eta_K(x)+ u\nabla_x\eta_K(x)\psi\ln{\frac{\psi}{\psi_\infty}})dxdR\\
\nonumber&=\int_{\mathbb{R}^{d}\times B}(-div (u\eta_K(x))\psi+ u\nabla_x\eta_K(x)\psi\ln{\frac{\psi}{\psi_\infty}})dxdR\\
\nonumber&=\int_{\mathbb{R}^{d}\times B}(-u \nabla_x\eta_K(x)\psi +u\nabla_x\eta_K(x)\psi\ln{\frac{\psi}{\psi_\infty}})dxdR\\
\nonumber&=\int_{\{K\leq |x|\leq 2K\} \times B}u \nabla_x\eta_K(x)\psi(\ln{\frac{\psi}{\psi_\infty}}-1)dxdR=J^K.
\end{align}
Using the fact that $\nabla_R\ln{\frac{\psi}{\psi_\infty}}=\nabla_R(\ln\psi-\ln{\psi_\infty})=\frac{\nabla_R\psi}{\psi}-\frac{\nabla_R\psi_\infty}{\psi_\infty}$, we deduce that
\begin{align}\label{3.18}
&\int_{\mathbb{R}^{d}\times B}\nabla{u}\cdot{R}\psi\nabla_{R}(\ln{\frac{\psi}{\psi_\infty}}\eta_K(x))dxdR\\
\nonumber&=\int_{\mathbb{R}^{d}\times B}\nabla{u}\cdot{R}\nabla_{R}\psi \eta_K(x) dxdR-\int_{\mathbb{R}^{d}\times B}\nabla{u}\cdot{R}\frac{\nabla_{R}\psi_\infty}{\psi_\infty}\eta_K(x)\psi dxdR\\
\nonumber&=\int_{\mathbb{R}^{d}\times B}-div_R(\nabla{u}\cdot{R})\psi\eta_K(x) dxdR-\int_{\mathbb{R}^{d}\times B}\nabla{u}\cdot{R}\frac{\nabla_{R}\psi_\infty}{\psi_\infty}\eta_K(x)\psi dxdR\\
\nonumber&=\int_{\mathbb{R}^{d}\times B}-(div~u)\psi\eta_K(x) dxdR-\int_{\mathbb{R}^{d}\times B}\nabla{u}\cdot{R}\frac{\nabla_{R}\psi_\infty}{\psi_\infty}\eta_K(x)\psi dxdR\\
\nonumber&=-\int_{\mathbb{R}^{d}\times B}\nabla{u}\cdot{R}\frac{\nabla_{R}\psi_\infty}{\psi_\infty}\eta_K(x)\psi dxdR.
\end{align}
Plugging (\ref{3.17}) and (\ref{3.18}) into (\ref{3.16}) yields
\begin{align}
4\int_{\mathbb{R}^{d}\times B}\psi_\infty\bigg|\nabla_{R}\sqrt{\frac{\psi}{\psi_{\infty}}}\bigg|^2\eta_K(x)dxdR=J^K-\int_{\mathbb{R}^{d}\times B}\nabla{u}\cdot{R}\frac{\nabla_{R}\psi_\infty}{\psi_\infty}\eta_K(x)\psi dxdR.
\end{align}
By virtue of the entropy condition, we deduce that
\begin{align}
J^K\leq \int_{\{K\leq |x|\leq 2K\} \times B}u \nabla_x\eta_K(x)\psi(\ln{\frac{\psi}{\psi_\infty}}-1)dxdR\lm \frac{\|u\|_{L^\infty}}{K} \int_{\{K\leq |x|\leq 2K\} \times B}\psi(\ln{\frac{\psi}{\psi_\infty}}-1)dxdR,
\end{align}
which leads to $\lim_{K\rightarrow\infty}J^K=0$. Combining with (3.4) and (3.16), we obtain
\begin{multline}
\nu\int_{|x|\leq K}|\nabla u|^2dx+4\int_{\{|x|\leq K\}\times B}\psi_\infty\bigg|\nabla_{R}\sqrt{\frac{\psi}{\psi_{\infty}}}\bigg|^2dxdR\leq I^K_1+I^K_2+I^K_3+I^K_4+J^k\\
-\int_{\mathbb{R}^d}\tau:\nabla u \eta_K  dx-\int_{\mathbb{R}^{d}\times B}\nabla{u}\cdot{R}\frac{\nabla_{R}\psi_\infty}{\psi_\infty}\eta_K\psi dxdR.
\end{multline}
Using the Fubini theorem, we have
\begin{align}
\int_{\mathbb{R}^d}\tau:\nabla u \eta_K dx=\sum_{1\leq i,j\leq d}\int_{\mathbb{R}^d\times B}R_i\partial_{R_j}\mathcal{U}\psi\partial_iu_j\eta_K dxdR.
\end{align}
Since $\psi_\infty=\frac{e^{-\mathcal{U}}}{\int_Be^{-\mathcal{U}}dR}$, it follow that
\begin{align}
\int_{\mathbb{R}^{d}\times B}\nabla{u}\cdot{R}\frac{\nabla_{R}\psi_\infty}{\psi_\infty}\eta_K\psi dxdR=-\sum_{1\leq i,j\leq d}\int_{\mathbb{R}^d\times B}\partial_iu_jR_i\partial_{R_j}\mathcal{U}\eta_K\psi dxdR.
\end{align}
Plugging (3.19) and (3.20) into (3.18) yields
\begin{align}
\nu\int_{|x|\leq K}|\nabla u|^2dx+4\int_{\{|x|\leq K\}\times B}\psi_\infty\bigg|\nabla_{R}\sqrt{\frac{\psi}{\psi_{\infty}}}\bigg|^2dxdR\leq I^K_1+I^K_2+I^K_3+I^K_4+J^k.
\end{align}
Passing the limit as $K$ goes to $\infty$, we deduce that
\begin{align}
\nu\int_{\mathbb{R}^d}|\nabla u|^2dx+4\int_{\mathbb{R}^d\times B}\psi_\infty\bigg|\nabla_{R}\sqrt{\frac{\psi}{\psi_{\infty}}}\bigg|^2dxdR=0,
\end{align}
which leads to $u=C$ and $\psi=f(x)\psi_\infty$ for some constant $C$ and function $f(x)$ respectively. Due to $\lim_{|x|\rightarrow\infty}u=0$, we obtain $u=0$. Moreover, $1=\int_B\psi dR=\int_Bf(x)\psi_\infty dR=f(x)\int_B\psi_\infty dR=f(x)$. Thus, we get $\psi=\psi_\infty$.
\subsection{The proof of Theorem \ref{th2}}

 Now we turn our attention to prove Theorem \ref{th2}. For $K>0$, choose $\eta^i_K(x_i)~(i=1,2,3)$ to be a positive smooth cut-off function satisfying:
\begin{align}
\eta^i_K(x_i)=1,~~\text{if}~~|x_i|\leq K,~~\eta^i_K(x_i)=0,~~\text{if}~~|x_i|\geq 2K,~~|\partial_i\eta^i_K(x_i)|\leq \frac{C}{K}, \text{for some constant}~C.
\end{align}
For each fixed $K>0$ we choose $v_K(x)=u(x)\eta^1_K(x_1)\eta^2_K(x_2)\eta^3_K(x_3)$ and $\phi_K(x)=\ln{\frac{\psi}{\psi_\infty}}\eta^1_K(x_1)\eta^2_K(x_2)\eta^3_K(x_3)$ as a test function with respective to $u$ and $\psi$, then we have
\begin{align}
\int_{\mathbb{R}^3}[(u\otimes u):\nabla (u\eta^1_K\eta^2_K\eta^3_K)]+P div(u\eta^1_K\eta^2_K\eta^3_K) dx =\int_{\mathbb{R}^3} \tau :\nabla (u \eta^1_K\eta^2_K\eta^3_K)+\nu\nabla u:\nabla (u\eta^1_K\eta^2_K\eta^3_K) dx,
\end{align}
\begin{align}
\int_{\mathbb{R}^{3}\times B}u\psi\cdot\nabla_{x}(\ln{\frac{\psi}{\psi_\infty}}\eta^1_K\eta^2_K\eta^3_K) dxdR =\int_{\mathbb{R}^{3}\times B}[-\nabla{u}\cdot{R}\psi+\psi_{\infty}\nabla_{R}\frac{\psi}{\psi_{\infty}}]\cdot\nabla_{R}(\ln{\frac{\psi}{\psi_\infty}}\eta^1_K\eta^2_K\eta^3_K) dxdR.
\end{align}
By a similar argument as in the proof of Theorem \ref{th1}, we deduce that
\begin{multline}\label{3.29}
\nu\int^K_{-K}\int^K_{-K}\int^K_{-K}|\nabla u|^2dx_1dx_2dx_3+4\int^K_{-K}\int^K_{-K}\int^K_{-K}\int_B\psi_\infty\bigg|\nabla_{R}\sqrt{\frac{\psi}{\psi_{\infty}}}\bigg|^2dx_1dx_2dx_3dR\\
\leq \overline{I}^K_1+\overline{I}^K_2+\overline{I}^K_3+\overline{I}^K_4+\overline{J}^K,
\end{multline}
where
\begin{align}
&\overline{I}^K_1=\nu\int_{\mathbb{R}^3}|\nabla u||u||\nabla (\eta^1_K\eta^2_K\eta^3_K)|dx, ~~\overline{I}^K_2=\frac{1}{2}\int_{\mathbb{R}^3}|u|^3||\nabla (\eta^1_K\eta^2_K\eta^3_K)|dx, \\ &\overline{I}^K_3=\int_{\mathbb{R}^3}|\tau||u||\nabla(\eta^1_K\eta^2_K\eta^3_K)|dx,~~ \overline{I}^K_4=\int_{\mathbb{R}^3}|P||u||\nabla (\eta^1_K\eta^2_K\eta^3_K)|dx, \\
&\overline{J}^K=\int_{\mathbb{R}^3 \times B}u \nabla_x(\eta^1_K\eta^2_K\eta^3_K)\psi(\ln{\frac{\psi}{\psi_\infty}}-1)dxdR.
\end{align}
Indeed, by the same token as the estimates for $I^K_1$, $I^K_3$ and $J^K$, one can obtain
\begin{align}
\lim_{K\rightarrow\infty}\overline{I}^K_1=0,~~ \lim_{K\rightarrow\infty}\overline{I}^K_3=0~~ \text{and}~~  \lim_{K\rightarrow\infty}\overline{J}^K=0.
\end{align}
Now we estimate $\overline{I}^K_2$ as follow.
\begin{align}
\overline{I}^K_2&\lm \int^{2K}_K\int^{2K}_{-2K}\int^{2K}_{-2K}|u|^3(\partial_1\eta^1_K)\eta^2_K\eta^3_Kdx_1dx_2dx_3
+\int^{-K}_{-2K}\int^{2K}_{-2K}\int^{2K}_{-2K}|u|^3(\partial_1\eta^1_K)\eta^2_K\eta^3_Kdx_1dx_2dx_3\\
\nonumber&
+\int^{2K}_{-2K}\int^{2K}_{K}\int^{2K}_{-2K}|u|^3\eta^1_K(\partial_2\eta^2_K)\eta^3_Kdx_1dx_2dx_3
+\int^{2K}_{-2K}\int^{-K}_{-2K}\int^{2K}_{-2K}|u|^3\eta^1_K(\partial_2\eta^2_K)\eta^3_Kdx_1dx_2dx_3\\
\nonumber&
+\int^{2K}_{-2K}\int^{2K}_{-2K}\int^{2K}_{K}|u|^3\eta^1_K\eta^2_K(\partial_3\eta^3_K)dx_1dx_2dx_3
+\int^{2K}_{-2K}\int^{2K}_{-2K}\int^{-K}_{-2K}|u|^3\eta^1_K\eta^2_K(\partial_3\eta^3_K)dx_1dx_2dx_3\\
\nonumber&=\overline{I}^K_{21}+\overline{I}^K_{22}+\overline{I}^K_{23}+\overline{I}^K_{24}+\overline{I}^K_{25}+\overline{I}^K_{26}.
\end{align}
We only treat with the term $\overline{I}^K_{21}$, and the others term can be estimated by the similar way. By virtue of H\"{o}lder's inequality, we get
\begin{align}
\overline{I}^K_{21}&\lm \frac{1}{K} \int^{2K}_K \int^{2K}_{-2K}\int^{2K}_{-2K}|u|^3dx_1dx_2dx_3\lm \sum_{1\leq i\leq 3}\frac{1}{K}\int^{2K}_K \int^{2K}_{-2K} \bigg(\int^{2K}_{-2K}|u_i|^{r_i}dx_3\bigg)^{\frac{3}{r_i}} K^{1-\frac{3}{r_i}}dx_1dx_2 \\
\nonumber&\lm \sum_{1\leq i\leq 3}\frac{ K^{1-\frac{3}{r_i}}K^{1-\frac{3}{q_i}}}{K}\int^{2K}_K \bigg[\int^{2K}_{-2K} \bigg(\int^{2K}_{-2K}|u_i|^{r_i}dx_3\bigg)^{\frac{q_i}{r_i}}dx_2\bigg]^{\frac{3}{q_i}} dx_1\\
\nonumber&\lm \sum_{1\leq i\leq 3} \frac{ K^{1-\frac{3}{r_i}}K^{1-\frac{3}{q_i}}K^{1-\frac{3}{p_i}}}{K}\bigg\{\int^{2K}_K\bigg[\int^{2K}_{-2K}\bigg(\int^{2K}_{-2K}|u_i|^{r_i}dx_3\bigg)^{\frac{q_i}{r_i}}dx_2\bigg]^{\frac{p_i}{q_i}}dx_1\bigg\}^{\frac{3}{p_i}}.
\end{align}
Since $\frac{1}{p_i}+\frac{1}{q_i}+\frac{1}{r_i}=\frac{2}{3}$, it follows that $\frac{K^{1-\frac{3}{r_i}}K^{1-\frac{3}{q_i}}K^{1-\frac{3}{p_i}}}{K}=1$. Then we have
\begin{align}
\overline{I}^K_{21}\lm \sum_{1\leq i\leq 3}\int^{2K}_K \bigg\{\int^{2K}_K\bigg[\int^{2K}_{-2K}\bigg(\int^{2K}_{-2K}|u_i|^{r_i}dx_3\bigg)^{\frac{q_i}{r_i}}dx_2\bigg]^{\frac{p_i}{q_i}}dx_1\bigg\}^{\frac{3}{p_i}}.
\end{align}
Thanks to $u_i\in L^{p_i}_{x_1}L^{q_i}_{x_2}L^{r_i}_{x_3}~(i=1,2,3)$, we deduce that the right hand side of the above inequality goes to $0$ as $K\rightarrow\infty$. Hence, we verify that $\lim_{K\rightarrow\infty}\overline{I}^K_{21}=0$. Moreover, we can prove that $\lim_{K\rightarrow\infty}\overline{I}^K_{2i}=0~(1\leq i\leq 6)$ by the same token, and then $\lim_{K\rightarrow\infty}\overline{I}^K_{2}=0$. Using the fact that $P=\sum_{1\leq i,j\leq 3}\mathfrak{R}^i\mathfrak{R}^j(u_iu_j-\tau_{ij})$ and by the similar argument as the estimate for ${I}^K_{4}$, we infer that $\lim_{K\rightarrow\infty}\overline{I}^K_{4}=0$. Passing the limit as $K\rightarrow\infty$ in the both sides of (\ref{3.29}) yields that
\begin{align}
\nu\int_{\mathbb{R}^3}|\nabla u|^2dx+4\int_{\mathbb{R}^3\times B}\psi_\infty\bigg|\nabla_{R}\sqrt{\frac{\psi}{\psi_{\infty}}}\bigg|^2dxdR=0.
\end{align}
By a similar argument as in the proof of Theorem \ref{th1}, we see that $u=0$ and $\psi=\psi_{\infty}$.

\section{$L^2$ decay rate}
This section is devoted to study the long time behaviour for the velocity of the FENE dumbbell model. More precisely, we prove the $L^2$ decay for the solutions of the FENE dumbbell model and obtain the $L^2$ decay rate. Without loss of generality, we take $\nu=1$ throughout this section.

\subsection{Co-rotation case}
Firstly, we consider the co-rotation FENE dumbbell model, that is, $\sigma(u)=\nabla u-(\nabla u)^{T}$. The existence of the solutions in $L^2$ was established in \cite{Masmoudi.G,Schonbek2}. Then our main result can be stated as follows.
\begin{theo}\label{th3}
Let $(u,\psi)$ be a weak solution of (1.2) with the initial data $u_0\in L^2\cap L^1$ and $\psi_0$ satisfies $\psi_0-\psi_\infty\in L^2_x(\mathcal{L}^2)$ and $\int_B\psi_0=1$ $a.e.$ in $x$. Then there exists a constant $C$ such that
\begin{align}
\int_{\mathbb{R}^d\times B}\frac{|\psi-\psi_\infty|^2}{\psi_\infty}dxdR\leq C\exp{(-Ct)},
\end{align}
\begin{align}
\|u\|_{L^2}\leq C(1+t)^{-\frac{d}{4}},\quad \text{if} \quad d\geq3, \quad \|u\|_{L^2}\leq C_l\ln^{-l}(e+t),\quad \text{if}\quad d=2,
\end{align}
where $l>0$ is arbitrarily integer and $C_l$ is a constant dependent on $l$.
\begin{proof}
By density argument, we only need to prove the estimate for the smooth solution. 
Since $\psi_\infty=\displaystyle\frac{(1-|R|^2)^k}{\int_{B}(1-|R|^2)^kdR}=\displaystyle\frac{(1-|R|^2)^k}{C_0}$, it follows that
\begin{multline}\label{4.3}
div_R([(\nabla u-(\nabla u)^T]R\psi_\infty)=\sum_{i,j}\partial_{R_i}[(\partial_iu^j-\partial_ju^i)R_j\psi_\infty]
\\
=\sum_{i,j}(\partial_iu^j-\partial_ju^i)\delta_{ij}\psi_\infty+\sum_{i,j}\frac{2k(\partial_iu^j-\partial_ju^i)R_jR_i(1-|R|^2)^{k-1}}{C_0}
=0.
\end{multline}
By virtue of the second equation of (1.2), we have
\begin{equation}
(\psi-\psi_\infty)_t+(u\cdot\nabla)(\psi-\psi_\infty)=div_{R}[-\sigma(u)\cdot{R}(\psi-\psi_\infty)+\psi_{\infty}\nabla_{R}\frac{\psi-\psi_\infty}{\psi_{\infty}}].
\end{equation}
Multiplying $\frac{\psi-\psi_\infty}{\psi_\infty}$ by both sides of the above equation and integrating over $B$ with $R$, we obtain
\begin{multline}\label{4.5}
\frac{1}{2}\frac{d}{dt}\int_B\frac{|\psi-\psi_\infty|^2}{\psi_\infty}+\frac{1}{2}u\cdot\nabla_x\int_B\frac{|\psi-\psi_\infty|^2}{\psi_\infty}+\int_B\psi_\infty|\nabla_R(\frac{\psi-\psi_\infty}{\psi_\infty})|^2
\\
=\int_B\sigma(u)R(\psi-\psi_\infty)\nabla_R(\frac{\psi-\psi_\infty}{\psi_\infty}).
\end{multline}
Thanks to integration by parts and (\ref{4.3}), we see that
\begin{multline}\label{4.6}
\int_B\sigma(u)R(\psi-\psi_\infty)\nabla_R(\frac{\psi-\psi_\infty}{\psi_\infty})=\int_B\sigma(u)R\psi_\infty[\frac{1}{2}\nabla_R(\frac{\psi-\psi_\infty}{\psi_\infty})^2]\\
=-\frac{1}{2}\int_Bdiv_R([(\nabla u-(\nabla u)^T]R\psi_\infty)(\frac{\psi-\psi_\infty}{\psi_\infty})^2=0.
\end{multline}
Plugging (\ref{4.6}) into (\ref{4.5}) and using the fact that $div~ u=0$, we deduce that
\begin{align}\label{4.7}
\frac{1}{2}\frac{d}{dt}\int_{\mathbb{R}^d\times B}\frac{|\psi-\psi_\infty|^2}{\psi_\infty}+\int_{\mathbb{R}^d\times B}\psi_\infty|\nabla_R(\frac{\psi-\psi_\infty}{\psi_\infty})|^2=0.
\end{align}
By virtue of the equation (1.2), we have $\int_B\psi dR=\int_B\psi_0 dR=1$, which leads to $\int_B(\psi-\psi_\infty) dR=0$. Taking advantage of Lemma \ref{Lemma2}, we infer that
\begin{align}
\frac{1}{2}\frac{d}{dt}\int_{\mathbb{R}^d\times B}\frac{|\psi-\psi_\infty|^2}{\psi_\infty}+C\int_{\mathbb{R}^d\times B}\frac{|\psi-\psi_\infty|^2}{\psi_\infty}\leq 0,
\end{align}
which leads to
\begin{align}
\frac{d}{dt}\bigg[\exp{(Ct)}\int_{\mathbb{R}^d\times B}\frac{|\psi-\psi_\infty|^2}{\psi_\infty}\bigg]\leq 0\Rightarrow \int_{\mathbb{R}^d\times B}\frac{|\psi-\psi_\infty|^2}{\psi_\infty}\leq \exp{(-Ct)}\int_{\mathbb{R}^d\times B}\frac{|\psi_0-\psi_\infty|^2}{\psi_\infty}.
\end{align}
Since $\partial_x \psi_\infty=0$, it follows that $div \tau= div \int_{B}(R\otimes\nabla_{R}\mathcal{U})\psi dR=div \int_{B}(R\otimes\nabla_{R}\mathcal{U})(\psi-\psi_\infty) dR$. Then, we may assume that $\tau =\int_{B}(R\otimes\nabla_{R}\mathcal{U})(\psi-\psi_\infty) dR$ .
By the standard energy estimate for the Navier-Stokes equations, we get
\begin{align}
\frac{1}{2}\frac{d}{dt}\|u\|^2_{L^2}+\|\nabla u\|^2_{L^2}=-\int_{\mathbb{R}^d}\tau:\nabla u\leq \frac{1}{2}\|\nabla u\|_{L^2}+\frac{1}{2}\|\tau\|^2_{L^2}.
\end{align}
Using Lemmas \ref{Lemma2}-\ref{Lemma3}, we verify that
\begin{align}
\frac{d}{dt}\|u\|^2_{L^2}+\|\nabla u\|^2_{L^2}\leq\|\tau\|^2_{L^2}\leq C\int_{\mathbb{R}^d\times B}\psi_\infty|\nabla_R(\frac{\psi-\psi_\infty}{\psi_\infty})|^2.
\end{align}
Let $\lambda\geq 2C$ be a sufficient large constant. From the above inequality and (\ref{4.7}), we deduce that
\begin{align}\label{4.12}
\frac{d}{dt}(\lambda\|\psi-\psi_\infty\|^2_{\mathcal{L}^2}+\|u\|^2_{L^2})+\lambda\int_{\mathbb{R}^d\times B}\psi_\infty|\nabla_R(\frac{\psi-\psi_\infty}{\psi_\infty})|^2+\|\nabla u\|^2_{L^2}\leq 0.
\end{align}
Taking $\lambda=2C$, we have
\begin{align}
\|u\|^2_{L^2}\leq \|u_0\|^2_{L^2}+2C\|\psi_0-\psi_\infty\|^2_{\mathcal{L}^2}<\infty.
\end{align}
Assume that $f$ is a positive continuous function and $f'(t)>0$. From (\ref{4.12}), we have
\begin{multline}
\frac{d}{dt}(f(t)\lambda\|\psi-\psi_\infty\|^2_{\mathcal{L}^2}+f(t)\|\widehat{u}\|^2_{L^2})+\lambda f(t)\int_{\mathbb{R}^d\times B}\psi_\infty|\nabla_R(\frac{\psi-\psi_\infty}{\psi_\infty})|^2+f(t)\int_{\mathbb{R}^d}|\xi|^2|\widehat{u}|^2d\xi\\
\leq f'(t)\lambda\|\psi-\psi_\infty\|^2_{\mathcal{L}^2}+f'(t)\|\widehat{u}\|^2_{L^2}.
\end{multline}
 Setting $S(t)=\{\xi:f(t)|\xi|^2\leq f'(t)\}$, then we obtain
\begin{multline}\label{4.15}
\frac{d}{dt}(f(t)\lambda\|\psi-\psi_\infty\|^2_{\mathcal{L}^2}+f(t)\|\widehat{u}\|^2_{L^2})+\lambda f(t)\int_{\mathbb{R}^d\times B}\psi_\infty|\nabla_R(\frac{\psi-\psi_\infty}{\psi_\infty})|^2\\
\leq f'(t)\lambda\|\psi-\psi_\infty\|^2_{\mathcal{L}^2}+f'(t)\int_{S(t)}|\widehat{u}|^2d\xi.
\end{multline}
By virtue of (1.2), we get
\begin{align}\label{4.16}
\widehat{u}=e^{-t|\xi|^2}\widehat{u_0}+\int^t_0e^{-(t-s)|\xi|^2}i\xi \mathcal{F}(\mathbb{P}(u\otimes u)+\mathbb{P}\tau) ds,
\end{align}
where $\mathbb{P}$ stands for Leray's project operator. Using the fact that $|\widehat{f}|\leq \|f\|_{L^1}$, we have
\begin{multline}
|\widehat{u}|\leq e^{-t|\xi|^2}|\widehat{u_0}|+|\xi|\int^t_0\|u\|^2_{L^2}ds+|\xi|t^{\frac{1}{2}}(\int^t_0|\widehat{\tau}|^2ds)^{\frac{1}{2}}\\
\leq
\|u_0\|_{L^1}+|\xi|t(\|u_0\|^2_{L^2}+C\|\psi_0-\psi_\infty\|^2_{L^2})+|\xi|t^{\frac{1}{2}}(\int^t_0|\widehat{\tau}|^2 ds)^{\frac{1}{2}},
\end{multline}
which leads to
\begin{align}\label{4.18}
\int_{S(t)}|\widehat{u}|^2d\xi&\lm \int_{S(t)}d\xi+t^2\int_{S(t)}|\xi|^2d\xi+t\int_{S(t)}|\xi|^2(\int^t_0|\widehat{\tau}|^2ds)d\xi\\
\nonumber&\lm \int^{\sqrt{\frac{f'(t)}{f(t)}}}_0r^{d-1}dr+t^2\int^{\sqrt{\frac{f'(t)}{f(t)}}}_0r^{d+1}dr+t\frac{f'(t)}{f(t)}\int^t_0\|\tau\|^2_{L^2}ds\\
\nonumber&\lm (\frac{f'(t)}{f(t)})^{\frac{d}{2}}+t^2(\frac{f'(t)}{f(t)})^{\frac{d}{2}+1}+t\frac{f'(t)}{f(t)}\int^t_0\int_{\mathbb{R}^d\times B}\psi_\infty|\nabla_R(\frac{\psi-\psi_\infty}{\psi_\infty})|^2ds.
\end{align}
Taking $f(t)=(1+t)^d$, then $f'(t)=d(1+t)^{d-1}$ and we have
\begin{align}\label{4.19}
\int_{S(t)}|\widehat{u}|^2d\xi\lm (1+t)^{-\frac{d}{2}+1}+\int^t_0\int_{\mathbb{R}^d\times B}\psi_\infty|\nabla_R(\frac{\psi-\psi_\infty}{\psi_\infty})|^2ds.
\end{align}
Plugging (\ref{4.19}) into (\ref{4.15}) and using the fact that $\|\psi-\psi_\infty\|_{\mathcal{L}^2}\lm \exp{(-Ct)}$ yield that
\begin{multline}
\frac{d}{dt}((1+t)^d\lambda\|\psi-\psi_\infty\|^2_{\mathcal{L}^2}+(1+t)^d\|\widehat{u}\|^2_{L^2})+\lambda (1+t)^d\int_{\mathbb{R}^d\times B}\psi_\infty|\nabla_R(\frac{\psi-\psi_\infty}{\psi_\infty})|^2\\
\leq C(1+t)^{\frac{d}{2}}+C(1+t)^{d-1}\int^t_0\int_{\mathbb{R}^d\times B}\psi_\infty|\nabla_R(\frac{\psi-\psi_\infty}{\psi_\infty})|^2ds,
\end{multline}
which implies that
\begin{multline}
((1+t)^d\lambda\|\psi-\psi_\infty\|^2_{\mathcal{L}^2}+(1+t)^d\|u\|^2_{L^2})+\lambda \int^t_0 (1+t')^d\int_{\mathbb{R}^d\times B}\psi_\infty|\nabla_R(\frac{\psi-\psi_\infty}{\psi_\infty})|^2\\
\lm 1+ \int^t_0(1+t')^{\frac{d}{2}}dt'+\int^t_0(1+t')^{d-1}\int^{t'}_0\int_{\mathbb{R}^d\times B}\psi_\infty|\nabla_R(\frac{\psi-\psi_\infty}{\psi_\infty})|^2 dsdt'\\
\lm (1+t)^{\frac{d}{2}+1}+\int^t_0(1+t')^{d}\int_{\mathbb{R}^d\times B}\psi_\infty|\nabla_R(\frac{\psi-\psi_\infty}{\psi_\infty})|^2dt'.
\end{multline}
By taking $\lambda$ sufficiently large, we obtain
\begin{align}
\|u\|^2_{L^2}\lm (1+t)^{-\frac{d}{2}+1}.
\end{align}
If $d\geq 3$, from (\ref{4.16}) we have
\begin{align}
|\widehat{u}|&\leq e^{-t|\xi|^2}|\widehat{u_0}|+|\xi|\int^t_0\|u\|^2_{L^2}ds+|\xi|t^{\frac{1}{2}}(\int^t_0|\widehat{\tau}|^2ds)^{\frac{1}{2}}\\
\nonumber&\leq
\|u_0\|_{L^1}+C|\xi|\int^t_0(1+s)^{1-\frac{d}{2}}ds+|\xi|t^{\frac{1}{2}}(\int^t_0|\widehat{\tau}|^2 ds)^{\frac{1}{2}}\\
\nonumber&\leq \|u_0\|_{L^1}+C|\xi|\int^t_0(1+s)^{-\frac{1}{2}}ds+|\xi|t^{\frac{1}{2}}(\int^t_0|\widehat{\tau}|^2 ds)^{\frac{1}{2}}\\
\nonumber&=\|u_0\|_{L^1}+C|\xi|\sqrt{1+t}+|\xi|t^{\frac{1}{2}}(\int^t_0|\widehat{\tau}|^2 ds)^{\frac{1}{2}},
\end{align}
which leads to
\begin{align}
\int_{S(t)}|\widehat{u}|^2d\xi&\lm \int_{S(t)}d\xi+(t+1)\int_{S(t)}|\xi|^2d\xi+t\int_{S(t)}|\xi|^2(\int^t_0|\widehat{\tau}|^2ds)d\xi\\
\nonumber&\lm \int^{\sqrt{\frac{f'(t)}{f(t)}}}_0r^{n-1}dr+(t+1)\int^{\sqrt{\frac{f'(t)}{f(t)}}}_0r^{n+1}dr+t\frac{f'(t)}{f(t)}\int^t_0\|\tau\|^2_{L^2}ds\\
\nonumber&\lm (\frac{f'(t)}{f(t)})^{\frac{d}{2}}+(t+1)(\frac{f'(t)}{f(t)})^{\frac{d}{2}+1}+t\frac{f'(t)}{f(t)}\int^t_0\int_{\mathbb{R}^d\times B}\psi_\infty|\nabla_R(\frac{\psi-\psi_\infty}{\psi_\infty})|^2ds.
\end{align}
Taking $f(t)=(1+t)^d$, we have
\begin{align}\label{4.25}
\int_{S(t)}|\widehat{u}|^2d\xi\lm (1+t)^{-\frac{d}{2}}+\int^t_0\int_{\mathbb{R}^d\times B}\psi_\infty|\nabla_R(\frac{\psi-\psi_\infty}{\psi_\infty})|^2ds.
\end{align}
Plugging (\ref{4.25}) into (\ref{4.15}) yields that
\begin{multline}
\frac{d}{dt}((1+t)^d\lambda\|\psi-\psi_\infty\|^2_{\mathcal{L}^2}+(1+t)^d\|\widehat{u}\|^2_{L^2})+\lambda (1+t)^d\int_{\mathbb{R}^d\times B}\psi_\infty|\nabla_R(\frac{\psi-\psi_\infty}{\psi_\infty})|^2\\
\leq C(1+t)^{\frac{d}{2}-1}+C(1+t)^{d-1}\int^t_0\int_{\mathbb{R}^d\times B}\psi_\infty|\nabla_R(\frac{\psi-\psi_\infty}{\psi_\infty})|^2ds,
\end{multline}
which implies that
\begin{multline}
((1+t)^d\lambda\|\psi-\psi_\infty\|^2_{\mathcal{L}^2}+(1+t)^d\|u\|^2_{L^2})+\int^t_0\lambda (1+t')^d\int_{\mathbb{R}^d\times B}\psi_\infty|\nabla_R(\frac{\psi-\psi_\infty}{\psi_\infty})|^2\\
\lm 1+ \int^t_0(1+t')^{\frac{d}{2}-1}dt'+\int^t_0(1+t')^{d-1}\int^{t'}_0\int_{\mathbb{R}^d\times B}\psi_\infty|\nabla_R(\frac{\psi-\psi_\infty}{\psi_\infty})|^2 dsdt'\\
\lm (1+t)^{\frac{d}{2}}+\int^t_0(1+t')^d\int_{\mathbb{R}^d\times B}\psi_\infty|\nabla_R(\frac{\psi-\psi_\infty}{\psi_\infty})|^2 dt'.
\end{multline}
By taking $\lambda$ sufficient large, we obtain that
\begin{align}
\|u\|^2_{L^2}\lm (1+t)^{-\frac{d}{2}}.
\end{align}
Now we turn our attention to the case $d=2$. From (\ref{4.18}), we see that
\begin{align}
\int_{S(t)}|\widehat{u}|^2d\xi\lm \frac{f'(t)}{f(t)}+t^2(\frac{f'(t)}{f(t)})^2+t\frac{f'(t)}{f(t)}\int^t_0\int_{\mathbb{R}^2\times B}\psi_\infty|\nabla_R(\frac{\psi-\psi_\infty}{\psi_\infty})|^2ds.
\end{align}
Taking $f(t)=\ln^3{(e+t)}$, then $f'(t)=\frac{3\ln^2(e+t)}{e+t}$. Thus, we have
\begin{align}\label{4.30}
\int_{S(t)}|\widehat{u}|^2d\xi\lm \frac{1}{\ln^2(e+t)}+\frac{1}{\ln(e+t)}\int^t_0\int_{\mathbb{R}^2\times B}\psi_\infty|\nabla_R(\frac{\psi-\psi_\infty}{\psi_\infty})|^2ds.
\end{align}
Plugging (\ref{4.30}) into (\ref{4.15}) yields
\begin{multline}
\frac{d}{dt}(\ln^3(e+t)\lambda\|\psi-\psi_\infty\|^2_{\mathcal{L}^2}+\ln^3(e+t)\|\widehat{u}\|^2_{L^2})+\lambda \ln^3(e+t)\int_{\mathbb{R}^d\times B}\psi_\infty|\nabla_R(\frac{\psi-\psi_\infty}{\psi_\infty})|^2\\
\leq \frac{C }{(e+t)}+\frac{C\ln(e+t) }{(e+t)}\int^t_0\int_{\mathbb{R}^d\times B}\psi_\infty|\nabla_R(\frac{\psi-\psi_\infty}{\psi_\infty})|^2ds,
\end{multline}
which implies that
\begin{multline}
(\ln^3(e+t)\lambda\|\psi-\psi_\infty\|^2_{\mathcal{L}^2}+\ln^3(e+t)\|\widehat{u}\|^2_{L^2})+\lambda\int^t_0 \ln^3(e+t')\int_{\mathbb{R}^d\times B}\psi_\infty|\nabla_R(\frac{\psi-\psi_\infty}{\psi_\infty})|^2dt'\\
\lm 1+ \int^t_0\frac{1}{(e+t')}dt'+\int^t_0\frac{\ln(e+t') }{(e+t')}\int^{t'}_0\int_{\mathbb{R}^d\times B}\psi_\infty|\nabla_R(\frac{\psi-\psi_\infty}{\psi_\infty})|^2dsdt'\\
\lm \ln(e+t)+\int^t_0\ln(e+t')\int_{\mathbb{R}^d\times B}\psi_\infty|\nabla_R(\frac{\psi-\psi_\infty}{\psi_\infty})|^2dt'.
\end{multline}
By taking $\lambda$ sufficiently large, we obtain
\begin{align}
\|u\|^2_{L^2}\lm \ln^{-2}(e+t)\quad \text{or} \quad \|u\|_{L^2}\lm \ln^{-1}(e+t).
\end{align}
By induction, we assume that $\|u\|_{L^2}\lm \ln^{-l}(e+t)$ for some $l\geq 1$. From (\ref{4.16}), we have
\begin{align}
|\widehat{u}|&\leq e^{-t|\xi|^2}|\widehat{u_0}|+|\xi|\int^t_0\|u\|^2_{L^2}ds+|\xi|t^{\frac{1}{2}}(\int^t_0|\widehat{\tau}|^2ds)^{\frac{1}{2}}\\
\nonumber&\leq
\|u_0\|_{L^1}+C_l|\xi|\int^t_0\ln^{-2l}(e+t)ds+|\xi|t^{\frac{1}{2}}(\int^t_0|\widehat{\tau}|^2 ds)^{\frac{1}{2}}\\
\nonumber&\leq \|u_0\|_{L^1}+C_l|\xi|(e+t)\ln^{-2l}(e+t)+|\xi|t^{\frac{1}{2}}(\int^t_0|\widehat{\tau}|^2 ds)^{\frac{1}{2}},
\end{align}
which leads to
\begin{align}
\int_{S(t)}|\widehat{u}|^2d\xi&\lm \int_{S(t)}d\xi+(t+e)\ln^{-2l}(e+t)\int_{S(t)}|\xi|^2d\xi+t\int_{S(t)}|\xi|^2(\int^t_0|\widehat{\tau}|^2ds)d\xi\\
\nonumber&\lm \int^{\sqrt{\frac{f'(t)}{f(t)}}}_0rdr+(t+e)\ln^{-2l}(e+t)\int^{\sqrt{\frac{f'(t)}{f(t)}}}_0r^2dr+t\frac{f'(t)}{f(t)}\int^t_0\|\tau\|^2_{L^2}ds\\
\nonumber&\lm (\frac{f'(t)}{f(t)})+(t+e)\ln^{-2l}(e+t)(\frac{f'(t)}{f(t)})^{2}+t\frac{f'(t)}{f(t)}\int^t_0\int_{\mathbb{R}^d\times B}\psi_\infty|\nabla_R(\frac{\psi-\psi_\infty}{\psi_\infty})|^2ds
\end{align}
Taking $f(t)=\ln^{2l+3}(e+t)$, then $f'(t)=(2l+3)\frac{\ln^{2l+2}(e+t)}{e+t}$. Thus, we get
\begin{align}\label{4.36}
\int_{S(t)}|\widehat{u}|^2d\xi\lm \frac{1}{\ln^{2l+2}(e+t)}+ \frac{1}{\ln(e+t)}\int^t_0\int_{\mathbb{R}^d\times B}\psi_\infty|\nabla_R(\frac{\psi-\psi_\infty}{\psi_\infty})|^2ds.
\end{align}
Plugging (\ref{4.36}) into (\ref{4.15}) yields
\begin{multline}
\frac{d}{dt}(\ln^{2l+3}(e+t)\lambda\|\psi-\psi_\infty\|^2_{\mathcal{L}^2}+\ln^{2l+3}(e+t)\|\widehat{u}\|^2_{L^2})+\lambda \ln^{2l+3}(e+t)\int_{\mathbb{R}^d\times B}\psi_\infty|\nabla_R(\frac{\psi-\psi_\infty}{\psi_\infty})|^2\\
\leq \frac{C }{(e+t)}+\frac{C\ln^{2l+1}(e+t) }{(e+t)}\int^t_0\int_{\mathbb{R}^d\times B}\psi_\infty|\nabla_R(\frac{\psi-\psi_\infty}{\psi_\infty})|^2ds,
\end{multline}
which implies
\begin{multline}
(\ln^{2l+3}(e+t)\lambda\|\psi-\psi_\infty\|^2_{\mathcal{L}^2}+\ln^{2l+3}(e+t)\|\widehat{u}\|^2_{L^2})+\lambda\int^t_0 \ln^{2l+3}(e+t')\int_{\mathbb{R}^d\times B}\psi_\infty|\nabla_R(\frac{\psi-\psi_\infty}{\psi_\infty})|^2dt'\\
\lm 1+ \int^t_0\frac{1}{(e+t')}dt'+\int^t_0\frac{\ln^{2l+1}(e+t') }{(e+t')}\int^{t'}_0\int_{\mathbb{R}^d\times B}\psi_\infty|\nabla_R(\frac{\psi-\psi_\infty}{\psi_\infty})|^2dsdt'\\
\lm \ln(e+t)+\int^t_0\ln^{2l+1}(e+t')\int_{\mathbb{R}^d\times B}\psi_\infty|\nabla_R(\frac{\psi-\psi_\infty}{\psi_\infty})|^2dt'.
\end{multline}
By taking $\lambda$ sufficiently large, we obtain
\begin{align}
\|u\|^2_{L^2}\lm \ln^{-2l-2}(e+t)\quad \text{or} \quad \|u\|_{L^2}\lm \ln^{-(1+1)}(e+t).
\end{align}
Therefore, by induction argument we already prove that for any $l\geq 1$
\begin{align}
\|u\|_{L^2}\leq C_l\ln^{-1}(e+t).
\end{align}
\end{proof}
\end{theo}
\begin{rema}
The above theorem improves the result obtained in \cite{Schonbek2} and we get the decay rate in dimension two.
\end{rema}
\subsection{General case}
Now we turn our attention to the general case, that is, $\sigma(u)=\nabla u$. For this propose, let us recall the global existence of strong solutions for (1.2) with small initial data.
\begin{theo}\cite{Masmoudi.G}\label{Masmoudi}
Let $s>1+\frac{d}{2}$. Assume $u_0\in H^s(\mathbb{R}^d)$ and $\psi_0-\psi_\infty\in H^s(\mathbb{R}^d;\mathcal{L}^2)$ with $\int_B\psi_0dR=1$ $a.e.$ in $x$. If there is a constant $\ep_0$ such that
\[\|u_0\|^2_{H^s}+\|\psi_0-\psi_\infty\|^2_{H^s(\mathcal{L}^2)}\leq \ep_0,\]
then there exists a unique global solution $(u,\psi)$ of (1.2) such that
 $u\in C(\mathbb{R}^+;H^s)\cap L^2_{loc}(\mathbb{R}^+;H^{s+1})$ and $\psi-\psi_\infty\in C(\mathbb{R}^+;H^s(\mathbb{R}^d;\mathcal{L}^2))\cap L^2_{loc}(\mathbb{R}^+;H^s(\mathbb{R}^d;\mathcal{H}^1))$. Moreover
 \[\|u\|^2_{H^s}+\|\psi-\psi_\infty\|^2_{H^s(\mathcal{L}^2)}\leq C\ep_0.\]
\end{theo}
Our result is stated as follows.
\begin{theo}
Assume that $(u,\psi)$ is the strong solution of (1.2) with the initial data $(u_0,\psi_0)$ under the condition of Theorem \ref{Masmoudi}. In addition, if $u_0\in L^1(\mathbb{R}^d)$ and $\sup_{R}\|\psi_0-\psi_\infty\|_{L^1}<\infty$, then there exists a constant $C$ such that
\[\begin{cases}
\displaystyle\int_{\mathbb{R}^d}|u|^2dx+\displaystyle\int_{\mathbb{R}^d\times B}\frac{|\psi-\psi_\infty|^2}{\psi_\infty}dxdR\leq C(1+t)^{-\frac{d}{2}+1} & \text{if}\quad d\geq 3,\\
\displaystyle\int_{\mathbb{R}^d}|u|^2dx+\displaystyle\int_{\mathbb{R}^d\times B}\frac{|\psi-\psi_\infty|^2}{\psi_\infty}dxdR\leq C\ln^{-1}(1+t)  & \text{if}\quad d=2.
\end{cases}
\]
\begin{proof}
By the standard $L^2$ energy method similar to Theorem \ref{th3}, we have
\begin{align}\label{4.41}
\frac{1}{2}\frac{d}{dt}\int_{\mathbb{R}^d}|u|^2dx+\int_{\mathbb{R}^d}|\nabla u|^2dx=-\int_{\mathbb{R}^d}\tau^{ij}\partial_iu^jdx=-\sum_{1\leq i,j\leq d}\int_{\mathbb{R}^d\times B}\partial_iu_jR_i\partial_{R_j}\mathcal{U}\psi dxdR,
\end{align}
\begin{multline}\label{4.42}
\frac{1}{2}\frac{d}{dt}\int_{\mathbb{R}^d\times B}\frac{|\psi-\psi_\infty|^2}{\psi_\infty}+\int_{\mathbb{R}^d\times B}\psi_\infty|\nabla_R(\frac{\psi-\psi_\infty}{\psi_\infty})|^2
\\
=\int_{\mathbb{R}^d\times B}\nabla uR(\psi-\psi_\infty)\nabla_R(\frac{\psi-\psi_\infty}{\psi_\infty})+\int_{\mathbb{R}^d\times B}\nabla u R\psi_\infty\nabla_R (\frac{\psi-\psi_\infty}{\psi_\infty}).
\end{multline}
By virtue of integration by parts and using the fact that $-\frac{\partial_{R_j}\psi_\infty}{\psi_\infty}=\partial_{R_j}\mathcal{U}$, we see that
\begin{align}
&\int_{\mathbb{R}^d\times B}\nabla u R\psi_\infty\nabla_R (\frac{\psi-\psi_\infty}{\psi_\infty})=\int_{\mathbb{R}^d\times B}\nabla u R\psi_\infty\nabla_R (\frac{\psi}{\psi_\infty})=-\int_{\mathbb{R}^d\times B}div_R(\nabla u R \psi_\infty)\frac{\psi}{\psi_\infty}\\
\nonumber&=-\int_{\mathbb{R}^d\times B} div u~\psi-\sum_{1\leq i,j\leq d}\int_{\mathbb{R}^d\times B} \partial_iu^jR_i(\partial_{R_j}\psi_\infty)\frac{\psi}{\psi_\infty}=\sum_{1\leq i,j\leq d}\int_{\mathbb{R}^d\times B}\partial_iu_jR_i\partial_{R_j}\mathcal{U}\psi.
\end{align}
Combining with (\ref{4.41}) and (\ref{4.42}) yields
\begin{multline}
\frac{1}{2}\frac{d}{dt}(\int_{\mathbb{R}^d}|u|^2dx+\int_{\mathbb{R}^d\times B}\frac{|\psi-\psi_\infty|^2}{\psi_\infty})+\int_{\mathbb{R}^d}|\nabla u|^2dx+\int_{\mathbb{R}^d\times B}\psi_\infty|\nabla_R(\frac{\psi-\psi_\infty}{\psi_\infty})|^2
\\
=\int_{\mathbb{R}^d\times B}\nabla uR(\psi-\psi_\infty)\nabla_R(\frac{\psi-\psi_\infty}{\psi_\infty}).
\end{multline}
Taking advantage of Cauchy-Schwarz's inequality and Lemma \ref{Lemma2}, we verify that
\begin{multline}
\frac{1}{2}\frac{d}{dt}(\int_{\mathbb{R}^d}|u|^2dx+\int_{\mathbb{R}^d\times B}\frac{|\psi-\psi_\infty|^2}{\psi_\infty})+\int_{\mathbb{R}^d}|\nabla u|^2dx+\int_{\mathbb{R}^d\times B}\psi_\infty|\nabla_R(\frac{\psi-\psi_\infty}{\psi_\infty})|^2
\\
\leq C\|\nabla u\|_{L^\infty}\int_{\mathbb{R}^d\times B}\psi_\infty|\nabla_R(\frac{\psi-\psi_\infty}{\psi_\infty})|^2.
\end{multline}
Using the fact that $\|\nabla u\|_{L^\infty}\leq \|u\|_{H^s}\leq C\ep_0$ with $\ep_0$ small enough, we obtain
\begin{align}
\frac{1}{2}\frac{d}{dt}(\int_{\mathbb{R}^d}|u|^2dx+\int_{\mathbb{R}^d\times B}\frac{|\psi-\psi_\infty|^2}{\psi_\infty})+\int_{\mathbb{R}^d}|\nabla u|^2dx+\frac{1}{2}\int_{\mathbb{R}^d\times B}\psi_\infty|\nabla_R(\frac{\psi-\psi_\infty}{\psi_\infty})|^2\leq 0.
\end{align}
Assume that $f$ is a positive continuous function and $f'(t) > 0$. From the above inequality, we see that
 \begin{multline}
\frac{d}{dt}(f(t)\|\psi-\psi_\infty\|^2_{\mathcal{L}^2}+f(t)\|\widehat{u}\|^2_{L^2})+ f(t)\int_{\mathbb{R}^d\times B}\psi_\infty|\nabla_R(\frac{\psi-\psi_\infty}{\psi_\infty})|^2+2f(t)\int_{\mathbb{R}^d}|\xi|^2|\widehat{u}|^2d\xi\\
\leq f'(t)\|\psi-\psi_\infty\|^2_{\mathcal{L}^2}+f'(t)\|\widehat{u}\|^2_{L^2}.
\end{multline}
 Setting $S(t)=\{\xi:2f(t)|\xi|^2\leq f'(t)\}$, then we obtain
\begin{multline}
\frac{d}{dt}(f(t)\|\psi-\psi_\infty\|^2_{\mathcal{L}^2}+f(t)\|\widehat{u}\|^2_{L^2})+ f(t)\int_{\mathbb{R}^d\times B}\psi_\infty|\nabla_R(\frac{\psi-\psi_\infty}{\psi_\infty})|^2\\
\leq f'(t)\|\psi-\psi_\infty\|^2_{\mathcal{L}^2}+f'(t)\int_{S(t)}|\widehat{u}|^2d\xi.
\end{multline}
Denote that $z=1-|R|$. By a simple calculation, we have $\int_B\frac{1}{z}dR=\int_S\int^1_0z^{d-2}dzd\sigma<\infty$ which implies that
\begin{align}
\widehat{\tau}\lm\int_B\frac{\mathcal{F}(\psi-\psi_\infty)}{1-|z|}\lm \sup_{R}|\mathcal{F}(\psi-\psi_\infty)|\lm \sup_{R}\|\psi-\psi_\infty\|_{L^1}\lm \sup_{R}\|\psi_0-\psi_\infty\|_{L^1},
\end{align}
together with (\ref{4.16}), we get
\begin{align}
|\widehat{u}|\leq e^{-t|\xi|^2}|\widehat{u_0}|+|\xi|\int^t_0\|u\|^2_{L^2}ds+|\xi|\int^t_0|\widehat{\tau}| ds
\lm
1+|\xi|t,
\end{align}
which leads to
\begin{align}
\int_{S(t)}|\widehat{u}|^2d\xi&\lm \int_{S(t)}d\xi+t^2\int_{S(t)}|\xi|^2d\xi\lm \int^{\sqrt{\frac{f'(t)}{f(t)}}}_0r^{d-1}dr+t^2\int^{\sqrt{\frac{f'(t)}{f(t)}}}_0r^{d+1}dr\\
\nonumber&\lm (\frac{f'(t)}{f(t)})^{\frac{d}{2}}+t^2(\frac{f'(t)}{f(t)})^{\frac{d}{2}+1}.
\end{align}
Taking $f(t) = (\eta+ t)^d$, where $\eta$ is a constant determinate later, then $f'(t)=d(\eta + t)^{d-1}$. Thus, we have
\begin{multline}
\frac{d}{dt}((\eta+ t)^d\|\psi-\psi_\infty\|^2_{\mathcal{L}^2}+(\eta+ t)^d\|\widehat{u}\|^2_{L^2})+ (\eta+ t)^d\int_{\mathbb{R}^d\times B}\psi_\infty|\nabla_R(\frac{\psi-\psi_\infty}{\psi_\infty})|^2\\
\leq d(\eta+t)^{d-1}\|\psi-\psi_\infty\|^2_{\mathcal{L}^2}+C(\eta + t)^{\frac{d}{2}}.
\end{multline}
By virtue of Lemma \ref{Lemma2} and taking $\eta$ large enough, we verify that
 \begin{align}
\frac{d}{dt}((\eta+ t)^d\|\psi-\psi_\infty\|^2_{\mathcal{L}^2}+(\eta+ t)^d\|\widehat{u}\|^2_{L^2})
\leq C(\eta + t)^{\frac{d}{2}},
\end{align}
which leads to
\begin{align}
\displaystyle\int_{\mathbb{R}^d}|u|^2dx+\displaystyle\int_{\mathbb{R}^d\times B}\frac{|\psi-\psi_\infty|^2}{\psi_\infty}dxdR\leq C(1+t)^{-\frac{d}{2}+1}.
\end{align}
If $d=2$, by taking $f(t)=\ln^3(\eta+t)$ and repeating the argument as above we infer that
\begin{align}
\displaystyle\int_{\mathbb{R}^d}|u|^2dx+\displaystyle\int_{\mathbb{R}^d\times B}\frac{|\psi-\psi_\infty|^2}{\psi_\infty}dxdR\leq C\ln^{-1}(1+t).
\end{align}
\end{proof}
\end{theo}
\begin{rema}
In the general case, one cannot obtain the $L^2$ energy estimate for the probability density, thus we cannot obtain the exponential decay rate. Moreover, the bootstrap argument as in the proof of Theorem \ref{th3} is invalid.
\end{rema}

\smallskip
\noindent\textbf{Acknowledgments} This work was partially supported by
NNSFC (No. 11271382 and No. 10971235), RFDP (No. 20120171110014),
and the key project of Sun Yat-sen University. 

\phantomsection
\addcontentsline{toc}{section}{\refname}

\end{document}